\documentclass[10pt]{amsart}
\usepackage[T1]{fontenc}
\usepackage{lmodern}
\usepackage[paperwidth=165mm,paperheight=240mm,textwidth=136mm,
  textheight=200mm,inner=16mm,top=25mm,headheight=12pt,headsep=5mm]{geometry}
\usepackage{amsmath,amsthm,amssymb,amscd,mathtools}
\usepackage{aliascnt}
\usepackage{hyperref,enumitem,cleveref}
\setlist[enumerate,1]{label=(\roman*),ref=(\roman*)}

\newtheorem{theorem}{Theorem}
\newaliascnt{lemma}{theorem}
\newtheorem{lemma}[lemma]{Lemma}
\aliascntresetthe{lemma}
\newaliascnt{proposition}{theorem}
\newtheorem{proposition}[proposition]{Proposition}
\aliascntresetthe{proposition}
\newaliascnt{corollary}{theorem}
\newtheorem{corollary}[corollary]{Corollary}
\aliascntresetthe{corollary}
\theoremstyle{definition}
\newaliascnt{definition}{theorem}
\newtheorem{definition}[definition]{Definition}
\aliascntresetthe{definition}
\newaliascnt{remark}{theorem}
\newtheorem{remark}[remark]{Remark}
\aliascntresetthe{remark}
\newaliascnt{example}{theorem}

\aliascntresetthe{example}

\usepackage{mathrsfs}
\usepackage[scr=boondoxo]{mathalpha}
\usepackage{xspace}
\hypersetup{
    colorlinks=true,
    linkcolor=purple,
    citecolor=magenta,
    filecolor=cyan,
    urlcolor=cyan,
    pdftitle={Nonassociative gradings on infinite-dimensional triangular algebras},
    }

\usepackage{xcolor}
\newcommand{\highlight}[1]{\textit{#1}}
\newcommand{\mathmode}[1]{\ensuremath{#1}\xspace}

\newcommand{\F}{\ensuremath{\mathbb{F}}\xspace} 
\newcommand{\NN}{\ensuremath{\mathbb{N}}\xspace} 

\newcommand{\G}{\ensuremath{G}\xspace} 
\renewcommand{\H}{\ensuremath{H}\xspace} 

\newcommand{\algbr}[1]{\ensuremath{\mathcal{#1}}\xspace} 
\newcommand{\aA}{\algbr{A}} 
\newcommand{\aB}{\algbr{B}}  
\newcommand{\aC}{\algbr{C}}  
\newcommand{\aD}{\algbr{D}}

\newcommand{\ideal}[1]{\ensuremath{\mathcal{#1}}\xspace} 
 
\newcommand{\iJ}{\ideal{J}}

\newcommand{\aset}[1]{\mathmode{\mathcal{#1}}} 
\newcommand{\sE}{\aset{E}}

\newcommand{\sI}{\aset{I}}

\newcommand{\sS}{\aset{S}}

\newcommand{\sY}{\aset{Y}}

\newcommand{\sset}{\mathmode{\star}-set\xspace} 

\newcommand{\B}{\ensuremath{\beta}\xspace} 

\newcommand{\V}{\ensuremath{\mathcal{V}}\xspace}

\newcommand{\range}{\operatorname{range}}
\newcommand{\End}{\operatorname{End}}

\newcommand{\EndF}{\End_{\F}}
\newcommand{\EndFV}{\ensuremath{\EndF\V}\xspace}

\newcommand{\id}{\operatorname{id}}

\newcommand{\UT}{\operatorname{UT}}
\newcommand{\UTB}{\UT_{\B}}

\newcommand{\UTBV}{\ensuremath{\UTB(\V)}\xspace}
\newcommand{\fin}{\operatorname{fin}}
\newcommand{\UTfin}{\UT^{\fin}}
\newcommand{\UTBlim}{\UT\limits_{\to\B}}
\newcommand{\UTBVlim}{\ensuremath{\UTBlim(\V)}\xspace}

\newcommand{\UTBfin}{\UTfin_{\B}}
\newcommand{\UTBVfin}{\ensuremath{\UTBfin(\V)}\xspace}
\newcommand{\UTnF}{\ensuremath{\UT_{n}(\F)}\xspace}
\newcommand{\UTIV}{\ensuremath{\UT_{I}(\V)}\xspace}

\newcommand{\Span}{\operatorname{Span}}
\newcommand{\Supp}{\operatorname{Supp}}

\DeclareMathOperator{\chr}{char}

\newcommand{\idemp}[1]{\mathmode{\mathscr{#1}}} 
\newcommand{\ee}{\idemp{e}}

\newcommand{\ff}{\idemp{f}}

\renewcommand{\gg}{\idemp{g}} 

\newcommand{\pp}{\idemp{p}}
\newcommand{\qq}{\idemp{q}}

\title[Nonassociative gradings on triangular algebras]{Nonassociative gradings on infinite-dimensional triangular algebras}

\author[P. Koshlukov]{Plamen Koshlukov}
\address{Department of Mathematics, State University of Campinas, Campinas, SP, Brazil}
\thanks{P.~Koshlukov was partially supported by FAPESP grant No.~2024/14914-9 and by CNPq grant No.~307184/2023-4}
\email{plamen@unicamp.br}

\author[W. Sch\"utzer]{Waldeck Sch\"utzer}
\address{Department of Mathematics, Universidade Federal de S\~{a}o Carlos, Brazil}
\email{waldeck@dm.ufscar.br}

\author[F. Yasumura]{Felipe Yukihide Yasumura}
\address{Department of Mathematics, State University of Campinas, Campinas, SP, Brazil}
\thanks{F.~Yasumura was supported by FAPESP grants No.~2024/14914-9 and No.~2025/16296-3.}
\email{yukihide@unicamp.br}


\begin{document}
\begin{abstract}
We study group gradings on triangular algebras relative to a
well-ordered basis, with multiplication $x\circ_\lambda y=xy+\lambda yx$.
For $\lambda\notin\{0,\pm1\}$, every grading on a finite-dimensional
triangular algebra is isomorphic to an elementary grading. In infinite
dimension, the same conclusion holds for the algebra of finite-range
maps, for the finite-map algebra on a basis of order type $\omega$,
and for closed gradings on the full triangular algebra. The
infinite-dimensional results also hold for $\lambda=1$ in characteristic
different from $2$.
For intermediate algebras whose finite-range ideal is graded, we obtain
an elementary description relative to a unitriangular basis.
The proofs construct neutral associative idempotents with singleton
diagonal supports. We also prove that the radical of the full triangular
algebra is graded for every nonzero $\lambda$, and that gradings on
the finite-range algebra and closed gradings on the full algebra have
finite support under the classification hypotheses.
\end{abstract}
\maketitle

\section{Introduction}
The classification of group gradings on simple algebras developed from
the study of Lie gradings \cite{PZ89} and their relation to associative
and Jordan gradings \cite{BaBrSh,BaKoch,BaSeZa2001}. The monograph
\cite{EK13} gives an account of this theory. Gradings on triangular
and incidence algebras provide a class of non-simple examples
\cite{BesDas,VinKoVa2004,KochY,KY2017,KY,SSY,VZ2007}.
The infinite-dimensional case requires different methods
\cite{BaBr,BaBrKoch,BaZa2010,KochY2,WY}.

Following \cite{WY}, we study triangular algebras relative to a
well-ordered basis, with multiplication $x\circ_\lambda y=xy+\lambda yx$.
We consider the full triangular algebra and its subalgebras of maps
with finite coordinate support or finite range. For
$\lambda\notin\{0,\pm1\}$, we first prove that every grading on a
finite-dimensional triangular algebra is isomorphic to an elementary
grading. For algebras containing the finite-range algebra as a graded
ideal, we obtain an elementary description relative to a unitriangular
basis. This gives an isomorphism to an elementary grading on the
original algebra in the finite-range and full triangular cases.
For the finite-map algebra on a basis of order type $\omega$, a
separate argument with neutral local units gives the same conclusion.
These infinite-dimensional results also hold for the Jordan product
when $\chr\F\neq2$. The final section treats closed gradings on
the full triangular algebra without assuming that the finite-range
ideal is graded. The Lie case is not treated.

The proofs reduce the classification to the construction of neutral
associative idempotents with singleton diagonal supports. Such
idempotents yield a triangular basis whose coordinate projections are
homogeneous. Finite-dimensional Jordan gradings may have primitive
neutral idempotents with two-point supports. We exclude this possibility
in infinite dimension by comparing corners isomorphic to $\UT_4(\F)$ and using
the well-ordering of the basis.

\section{Preliminaries}
\subsection{Gradings and graded algebras}
Let $G$ be a group and $\mathcal{A}$ a vector space over a field $\mathbb{F}$. A \emph{$G$-grading} on $\mathcal{A}$ is a vector space decomposition $\mathcal{A}=\bigoplus_{g\in G}\mathcal{A}_g$. We say that $\mathcal{A}$ is \emph{$G$-graded} if it is endowed with a $G$-grading. The vector subspace $\mathcal{A}_g$ is called the \emph{homogeneous component} of degree $g$, and its nonzero elements are called homogeneous of degree $g$. When the decomposition is not explicitly given, we shall write $(\mathcal{A})_g$ to denote the homogeneous component of degree $g$. Given $0\ne x\in\mathcal{A}_g$, we denote $\deg_Gx=g$, or simply $\deg x=g$ when there is no risk of ambiguity. The \emph{support} of a grading $\Gamma$ is defined as
\[
\mathrm{Supp}\,\Gamma=\{g\in G\mid\mathcal{A}_g\ne0\}.
\]
A vector subspace $\mathcal{S}\subseteq\mathcal{A}$ is said to be \emph{graded} (or \textsl{homogeneous}) if $\mathcal{S}=\bigoplus\mathcal{A}_g\cap\mathcal{S}$. A graded subspace naturally has an induced $G$-grading.

Given two $G$-graded vector spaces $\mathcal{A}=\bigoplus_{g\in G}\mathcal{A}_g$ and $\mathcal{C}=\bigoplus_{g\in G}\mathcal{C}_g$, a $G$-graded map is a linear map $\psi:\mathcal{A}\to\mathcal{C}$ such that $\psi(\mathcal{A}_g)\subseteq\mathcal{C}_g$, for all $g\in G$.

Now, assume that $\mathcal{A}$ is an $\Omega$-algebra, for a signature $\Omega$ containing multilinear operations. We say that $\mathcal{A}$ is a \emph{$G$-graded algebra} if $\mathcal{A}$ is endowed with a $G$-grading such that
$$
\omega(a_1,\ldots,a_m)\in\left(\mathcal{A}\right)_{g_1\cdots g_m},
$$
for each $a_i\in\left(\mathcal{A}\right)_{g_i}$, $i=1,\ldots,m$, and each $\omega=\omega(x_1,\ldots,x_m)\in\Omega$.

Two $G$-graded $\Omega$-algebras $\mathcal{A}$ and $\mathcal{B}$ are said to be isomorphic if there is a $G$-graded linear map $\mathcal{A}\to\mathcal{B}$ that is an isomorphism of $\Omega$-algebras. In this case, we denote $\mathcal{A}\cong\mathcal{B}$.

If $\Gamma$ and $\Gamma'$ are $G$-gradings on the $\Omega$-algebra $\mathcal{A}$, then we say that $\Gamma$ and $\Gamma'$ are isomorphic, denoted by $\Gamma\cong\Gamma'$ if $(\mathcal{A},\Gamma)\cong(\mathcal{A},\Gamma')$.

\subsection{Gradings on upper triangular matrix algebras}

Throughout this subsection, $n$ stands for a positive integer.
The group gradings on the associative algebra $\UT_n(\F)$ were
classified by Valenti and Zaicev in \cite{VZ2007}. They proved that
every group grading on $\UT_n(\F)$ is isomorphic to an elementary
grading, in which every matrix unit is homogeneous.

More precisely, if
\[
\gamma:\{1,\ldots,n\}\longrightarrow G,
\]
then assigning
\[
\deg(\ee_{ij})=\gamma(i)\gamma(j)^{-1},
\qquad 1\le i\le j\le n,
\]
defines an elementary $G$-grading on $\UT_n(\F)$.
Conversely, every group grading on $\UT_n(\F)$ is isomorphic to one of
this form.

\subsection{Topological gradings}
We shall use the following notions for gradings on topological algebras.

\begin{definition}[{\cite[Definition 9]{WY}}]
Let $\Gamma\colon\mathcal{A}=\bigoplus_{g\in G}\mathcal{A}_g$ be a $G$-grading on the topological algebra $\mathcal{A}$, and denote by $\pi_g\colon\mathcal{A}\to\mathcal{A}$ the respective projection of $\mathcal{A}$ over $\mathcal{A}_g$.
\begin{enumerate}
\item We say that $\Gamma$ is \emph{closed} if each subspace $\mathcal{A}_g$ is a (topologically) closed subspace.
\item We say that $\Gamma$ is \emph{continuous} if each $\pi_g$ is a continuous map.
\end{enumerate}
\end{definition}

\section{Triangular algebras and elementary gradings}

The zero algebra has only the trivial grading. Henceforth, let \V be a
nonzero vector space over \F with a well-ordered basis \B. The
\highlight{finite topology} on $\EndFV$ has as a neighborhood basis of
zero the subspaces
\[
U(W)=\{T\in\EndFV\mid T(W)=0\},
\]
where $W$ ranges over the finite-dimensional subspaces of \V. All
subalgebras of $\EndFV$ are endowed with the induced topology.

For each $v\in\B$, let $\V_{\leq v}= \Span\{u\in\B\mid u\leq v\}$. A linear map $T\in\EndFV$ is \highlight{triangular} (with respect to \B) if $T(v)\in \V_{\leq v}$, for every $v\in\B$. We denote by
\[\UTBV=\{T\in\EndFV\mid T(v)\in\V_{\leq v},\forall v\in\B\}\] the subset of \EndFV of all maps that are triangular with respect to \B. It is a unital subalgebra of \EndFV, closed in the finite topology, and is a \highlight{triangularizable} algebra in the sense of Mesyan \cite[Definition~2.1]{Mesb}.

For every pair $u\leq v$ in \B, let $\ee_{uv}\in\EndFV$ be the linear map defined on the basis by
\[
\ee_{uv}(w)=\begin{cases}
u,& \text{if }w=v,\\
0,& \text{otherwise}.
\end{cases}
\]
These maps are called the \highlight{matrix units} relative to~\B. In particular, $\ee_{vv}$ is a projection of \V onto the one-dimensional subspace $\F v$. 

For $T\in\EndFV$, let $T(u,v)$ denote the coefficient of $u$ in
$T(v)$. For $u\leq v$, we have
\[
\ee_{uu}T\ee_{vv}=T(u,v)\ee_{uv}.
\]
Every $T\in\UTBV$ is uniquely determined by its coordinates. More precisely,
\[
T=\sum_{u\leq v}T(u,v)\ee_{uv},
\]
where the sum converges in the finite topology. Thus, the matrix units form a \highlight{topological basis} of $\UTBV$.

The \highlight{\B-support} of a linear map $T\in\EndFV$ is the set
\[
\Supp_\B(T)=\{v\in\B\mid T(v,v)\neq0\}.
\]
When the basis $\B$ is clear from the context, we omit the subscript $\B$, and we simply write $\Supp(T)$ and refer to it as the \highlight{support} of $T$.

When $\B$ is countable, we shall often enumerate its elements as
\[
\B=\{v_i\mid i\in\NN\}.
\]
In this case, we abbreviate $\ee_{v_i,v_j}$ by $\ee_{ij}$ and write $T(i,j)$ in place of $T(v_i,v_j)$.

The notation for matrix units and coordinates depends on the choice of the basis~$\B$. Since the basis will always be specified or clear from the context, we shall suppress this dependence throughout the paper.

For each finite subset $I\subseteq\B$, we define
\[
\UTIV=\Span\{\ee_{uv}\mid u,v\in I,\ u\leq v\}.
\]
Then $\UTIV$ is a finite-dimensional subalgebra of $\UTBV$. Moreover, if
$n=|I|$, then
\[
\UTIV\cong\UTnF.
\]
We also define
\[
\UTBVlim=\Span\{\ee_{uv}\mid u,v\in\B,\ u\leq v\}.
\]
Its elements are precisely the \highlight{finite maps} in $\UTBV$, that is, the
maps having only finitely many nonzero coordinates. Thus $\UTBVlim$ is a
subalgebra of $\UTBV$. If $\B$ is infinite, then $\UTBVlim$ has no identity
element.

Every finite map belongs to $\UTIV$ for some finite subset $I\subseteq\B$.
Thus $\UTBVlim$ is the directed union of these finite-dimensional subalgebras.

We define $\UTBVfin$ to be the subalgebra of $\UTBV$ consisting of all maps $x$
such that $\dim_\F\range(x)<+\infty$.
Its elements are said to have \highlight{finite range}, and
\[
\UTBVlim\subseteq\UTBVfin.
\]
If $\B$ is infinite, then $\UTBVfin$ has no identity element.

Finally, if $\B$ is finite, then
\[
\UTBV=\UTBVfin=\UTBVlim\cong\UTnF,
\]
where $n=|\B|$.

Throughout the paper, uppercase letters $T,S,\ldots$ will denote elements of
$\EndFV$, while lowercase letters $x,y,\ldots$ will denote elements of
$\UTBV$, $\UTBVfin$, and $\UTBVlim$.

We shall study group gradings on $\UTBV$, $\UTBVlim$ and $\UTBVfin$
viewed as non-associative algebras, including the Jordan case. The Lie
case is not treated here. To this end, it is convenient to consider a
slightly more general setting.

Let $X=\{x_1,x_2,\ldots\}$ be a countable set of noncommuting variables, and let
$\F\langle X\rangle$ denote the free associative algebra on $X$.
We fix a set $\Omega\subseteq\F\langle X\rangle$ of multilinear polynomials in $\F\langle X\rangle$. If $\mathcal{A}$ is an associative algebra, then each $\omega\in\Omega$ induces a multilinear operation on $\mathcal{A}$ by evaluation. Thus, every associative algebra naturally becomes an $\Omega$-algebra.
We use the definition of a graded $\Omega$-algebra given in the
preliminaries for these induced operations.

For any map $\gamma\colon\B\to G$, we define a vector space $G$-grading
$\Gamma(\gamma)$ on $\UTBVlim$ by declaring each matrix unit $\ee_{uv}$ to be
homogeneous of degree
\[
\deg(\ee_{uv})=\gamma(u)\gamma(v)^{-1}.
\]
The grading $\Gamma(\gamma)$ is called \highlight{elementary}. For each
finite subset $I\subseteq\B$, the subalgebra $\UTIV$ is graded with respect to
$\Gamma(\gamma)$. We denote the induced grading on $\UTIV$ by
$\Gamma_I(\gamma)$.

More generally, let $\aA\subseteq\UTBV$ be an associative subalgebra and
let $\Gamma\colon\aA=\bigoplus_{g\in G}\aA_g$ be a $G$-grading. We say
that $\Gamma$ is \highlight{elementary relative to $\B$} if there exists
a map $\gamma\colon\B\to G$ such that
\[
\aA_g=
\left\{
x\in\aA\ \middle|\
x(u,v)\neq0
\Longrightarrow
\gamma(u)\gamma(v)^{-1}=g
\right\}
\]
for every $g\in G$. The map $\gamma$ is called a \highlight{defining
map} of $\Gamma$. For $\aA=\UTBVlim$, the grading $\Gamma(\gamma)$ is
elementary relative to $\B$ with defining map $\gamma$. If $\B'$ is a
well-ordered basis satisfying $\UTBV=\UT_{\B'}(\V)$, elementary relative
to $\B'$ is defined by the same condition using the coordinates relative
to $\B'$.

It is known that every elementary grading $\Gamma(\gamma)$ on
$\UTBVlim$ is continuous (see~\cite[Proposition~15]{WY}). Hence, if
\[
\Supp\Gamma(\gamma)
=\{g\in G\mid (\UTBVlim)_g\neq0\}
\]
is finite, then the closures of its homogeneous components form a
$G$-grading on $\UTBV$ (see \cite[Lemma 11]{WY}). The closure of
$(\UTBVlim)_g$ consists of the maps whose nonzero coordinates $x(u,v)$
all satisfy $\gamma(u)\gamma(v)^{-1}=g$, so this grading is elementary
relative to $\B$ with defining map $\gamma$.

Since $\UTBVlim$ is the union of the graded subalgebras $\UTIV$,
its elementary grading is an $\Omega$-grading if and only if every
finite restriction is an $\Omega$-grading. Indeed, any finite family of
finite maps lies in one such subalgebra. Moreover,
$\Supp\Gamma(\gamma)=\bigcup_{I\subseteq\B,\,|I|<\infty}
\Supp\Gamma_I(\gamma)$.

From now on, we specialize to
$\Omega_\lambda=\{xy+\lambda yx\}$, where $\lambda\in\F^\times$,
and write $x\circ_\lambda y=xy+\lambda yx$.
Powers for this product are defined by $x^{\circ_\lambda1}=x$ and
$x^{\circ_\lambda(k+1)}=x^{\circ_\lambda k}\circ_\lambda x$.
The case $\lambda=0$ is associative and was treated in \cite{WY}.
The values $\lambda=1$ and $\lambda=-1$ give the Jordan and Lie
products, respectively. We write $\circ$ for $\circ_1$.

The finite-dimensional Jordan and Lie cases of the next criterion
were proved in \cite[Lemma~2]{KY2017} and \cite[Lemma~10]{KY}.

\begin{proposition}\label{prop:elementary-omega}
Let $\lambda\in\F^\times$. The elementary grading
$\Gamma(\gamma)$ on $\UTBVlim$ is an $\Omega_\lambda$-grading if and
only if $\Supp\Gamma(\gamma)$ generates an abelian subgroup of $G$.
\end{proposition}

\begin{proof}
Let $v_0$ be the least element of $\B$ and put
$h_v=\gamma(v_0)\gamma(v)^{-1}$.
Suppose that $\Gamma(\gamma)$ is an $\Omega_\lambda$-grading.
For $v_0<u<v$, the nonzero product
\[
\ee_{uv}\circ_\lambda\ee_{v_0u}=\lambda\ee_{v_0v}
\]
gives $(h_u^{-1}h_v)h_u=h_v$, so $h_uh_v=h_vh_u$.
Since $h_{v_0}=1$, the elements $h_v$ commute pairwise.
Every degree in the support is of the form $h_u^{-1}h_v$, proving
necessity.

Conversely, the degrees in an elementary grading multiply under associative
multiplication of matrix units. If the support generates an abelian
subgroup, both summands of
$\ee_{ab}\circ_\lambda\ee_{cd}$ have degree
$\deg(\ee_{ab})\deg(\ee_{cd})$ whenever they are nonzero.
Multilinearity proves the assertion.
\end{proof}

For $\lambda^2\ne1$, associative multiplication is recovered by
\begin{equation}\label{eq:recover:product}
xy=\frac{x\circ_\lambda y-\lambda(y\circ_\lambda x)}{1-\lambda^2}.
\end{equation}

\begin{lemma}\label{product}
Let $\lambda\in\F\setminus\{0,\pm1\}$ and let $\Gamma$ be an
$\Omega_\lambda$-grading on an associative algebra $\aA$.
If $x\in\aA_g$, $y\in\aA_h$, and $gh=hg$, then $xy\in\aA_{gh}$.
\end{lemma}
\begin{proof}
Both terms in the numerator of \cref{eq:recover:product} belong to
$\aA_{gh}$.
\end{proof}

\begin{lemma}\label{hom_is_hom}
Let $\lambda\in\F\setminus\{0,\pm1\}$ and let $\aA$, $\aB$ be
associative algebras. A linear map $\varphi\colon\aA\to\aB$ is an
$\Omega_\lambda$-homomorphism if and only if it is an associative
homomorphism.
\end{lemma}
\begin{proof}
One direction follows from \cref{eq:recover:product}, and the other
from $x\circ_\lambda y=xy+\lambda yx$.
\end{proof}

\begin{lemma}\label{lem:identity:neutral}
Let $\lambda\in\F\setminus\{-1\}$, and let $\aA$ be a unital
associative algebra endowed with an $\Omega_\lambda$-grading.
Then $1\in\aA_1$.
\end{lemma}
\begin{proof}
Write $1=\sum_h a_h$, where $a_h\in\aA_h$.
For homogeneous $x\in\aA_g$, comparison of degrees in
$1\circ_\lambda x=(1+\lambda)x$ gives
$a_h\circ_\lambda x=0$ whenever $h\ne1$.
By linearity this holds for every $x\in\aA$.
Taking $x=1$ gives $(1+\lambda)a_h=0$, so $1=a_1$.
\end{proof}

\begin{remark}\label{rem:grad:inheritance}
Let $\aA=\bigoplus_{g\in\G}\aA_g$ be an associative algebra endowed
with an $\Omega_\lambda$-grading, and let $\ee\in\aA_1$ be an
associative idempotent.

If $\lambda\in\F\setminus\{\pm1\}$, then the associative grading
rule for $\lambda=0$, and \Cref{product} otherwise, give
$\aA_1\aA_g\subseteq\aA_g$ and
$\aA_g\aA_1\subseteq\aA_g$ for every $g\in\G$. In particular,
$\aA_1$ is closed under the associative multiplication and
$\ee\aA_g\ee\subseteq\aA_g$.

If $\lambda=1$ and $\chr\F\neq2$, define
\[
U_\ee\colon\aA\longrightarrow\aA,
\qquad
U_\ee(x)=\frac12\bigl(\ee\circ(\ee\circ x)-\ee\circ x\bigr)
=\ee x\ee.
\]
Since $\ee$ is homogeneous of degree $1$, the map $U_\ee$ preserves
each homogeneous component. Thus, in either case, the corner
$\ee\aA\ee$ inherits the $\Omega_\lambda$-grading, with homogeneous
components $\ee\aA_g\ee$.
\end{remark}

A nonzero associative idempotent in a subspace is \highlight{primitive}
in the subspace if it cannot be written as the sum of two nonzero orthogonal associative idempotents in that subspace. In particular,
primitivity in $\aA_1$ always has this meaning, including in the
Jordan case.

\begin{lemma}\label{lem:orthogonalize}
Let $\aA\subseteq\UTBV$ be an associative subalgebra endowed with an
$\Omega_\lambda$-grading. Suppose that $\ee,\ff\in\aA_1$ are
associative idempotents and $\ff\ee=0$. Then $\ee-\ee\ff$ and
$\ff-\ee\ff$ belong to $\aA_1$ and are associative idempotents orthogonal to
$\ff$ and $\ee$, respectively. Their supports equal those of $\ee$
and $\ff$, respectively, and
$\range(\ee-\ee\ff)=\range(\ee)$.
\end{lemma}
\begin{proof}
Since $\ee\ff=\ee\circ_\lambda\ff$, the product belongs to $\aA_1$.
The relation $\ff\ee=0$ gives
$(\ee-\ee\ff)^2=\ee-\ee\ff$ and
$(\ff-\ee\ff)^2=\ff-\ee\ff$, as well as the asserted orthogonality.
Since diagonal coordinates multiply, $\ff\ee=0$ gives
$\ee(v,v)\ff(v,v)=0$ for every $v$, so $\ee\ff$ has zero diagonal and the
supports are unchanged.
Finally, $\ee-\ee\ff=\ee(1-\ff)$ and
$(\ee-\ee\ff)\ee=\ee$ give equality of the ranges.
\end{proof}

\begin{lemma}\label{lem:radicals:finite-algebras}
If $\aA=\UTIV$, $\UTBVlim$, $\UTBVfin$, or $\UTBV$, then
\[
J(\aA)=\{x\in\aA\mid x(v,v)=0\text{ for every }v\in\B\}.
\]
For $\UTIV$ and $\UTBVlim$, this is the span of the matrix units
$\ee_{uv}\in\aA$ with $u<v$.
Every element of $J(\UTBV)$ is locally nilpotent.
\end{lemma}

\begin{proof}
The diagonal map is a surjective associative homomorphism from $\aA$
onto $\F^I$ if $\aA=\UTIV$, and onto the direct sum
$\F^{(\B)}$ if $\aA=\UTBVlim$ or $\UTBVfin$, and onto
$\F^\B$ if $\aA=\UTBV$. For $\aA=\UTBVfin$, the
diagonal support of a finite-range triangular transformation is finite,
because its images at distinct points of the diagonal support are
linearly independent by triangularity.
These diagonal algebras have zero Jacobson radical. Hence $J(\aA)$ is
contained in the kernel of the diagonal map.

Every strictly triangular transformation $z$ is locally nilpotent.
Indeed, by induction over $\B$, the finitely many basis vectors
occurring in $z(v)$ are annihilated by powers of $z$. Thus a power
of $z$ annihilates $v$.

Conversely, let $x\in\aA$ have zero diagonal and let
$y\in\F1+\aA$, where $1=\id_\V$.
The transformation $z=xy$ is strictly triangular. If $\aA=\UTBV$,
the pointwise finite sum $\sum_{k\geq0}z^k$ belongs to $\aA$ and
is the inverse of $1-z$. In the other cases, $z$ has finite range.
Choose a finite $I\subseteq\B$ with $\range(z)\subseteq\Span I$.
The restriction of $z$ to $\Span I$ is strictly triangular, so
$z^{|I|+1}=0$. Hence $1-z$ is invertible in $\F1+\aA$.
The Jacobson criterion gives $x\in J(\aA)$.
\end{proof}

\begin{lemma}\label{Jgraded}
Let $\lambda\in\F^\times$, and let $\Gamma$ be an
$\Omega_\lambda$-grading on $\aA$, where
$\aA=\UTIV$, $\UTBVlim$, or $\UTBVfin$. Denote by
\[
(a,b,c)_\lambda
=
(a\circ_\lambda b)\circ_\lambda c
-a\circ_\lambda(b\circ_\lambda c)
\]
the associator of the product $\circ_\lambda$ (that is the $\Omega_\lambda$-associator). Then
\[
J(\aA)
=
\Span\{(a,b,c)_\lambda\mid a,b,c\in\aA\}.
\]
Consequently, $J(\aA)$ is a graded ideal.
\end{lemma} 
\begin{proof}
Let $\mathcal S$ denote the linear span of all
$\Omega_\lambda$-associators. Since the associator is a multilinear and $G$-homogeneous
$\Omega_\lambda$-polynomial, $\mathcal S$ is graded.

By \Cref{lem:radicals:finite-algebras}, the diagonal map
$x\mapsto (x(v,v))_{v\in\B}$ induces an injective associative
homomorphism $\aA/J(\aA)\to\F^\B$. Hence $\aA/J(\aA)$ is commutative and
associative. The product induced by $\circ_\lambda$ on this quotient
is therefore associative, so
$\mathcal S\subseteq J(\aA)$.

Suppose first that $\aA=\UTIV$ or $\aA=\UTBVlim$. For every $\ee_{uv}\in\aA$ with $u<v$,
a direct calculation gives
\[
\ee_{uv}
=
\left(
\lambda^{-1}\ee_{uv},
\ee_{vv},
\ee_{uu}
\right)_\lambda.
\]
Since the matrix units $\ee_{uv}$, $u<v$, span $J(\aA)$, it follows
that $J(\aA)\subseteq\mathcal S$.

Now let $\aA=\UTBVfin$ and let $x\in J(\aA)$. There exists a finite set
$I\subseteq\B$ such that
$\range(x)\subseteq\Span I$. Put
$\pp=\sum_{v\in I}\ee_{vv}$. Then $\pp x=x$ and
\[
x=x(1-\pp)+\pp x\pp.
\]
Since $\pp x\pp\in\UTIV$ has zero diagonal, the finite-dimensional
case gives $\pp x\pp\in\mathcal S$.

Moreover, $(x,\pp,\pp)_\lambda
= \lambda(2\pp x\pp-x\pp-\pp x) =
-\lambda x(1-\pp)$. 
Thus we obtain 
$x(1-\pp)=-\lambda^{-1}(x,\pp,\pp)_\lambda\in\mathcal S$.
Consequently, $x\in\mathcal S$.

Therefore $J(\aA)=\mathcal S$, and since $\mathcal S$ is graded,
so is $J(\aA)$.
\end{proof}

A basis $\sY=\{y_v\mid v\in\B\}$ is
\highlight{unitriangular with respect to $\B$} if
$y_v-v\in\V_{<v}=\Span\{u\in\B\mid u<v\}$ for every $v\in\B$.
We give $\sY$ the order induced by its indices in $\B$.

\begin{lemma}\label{lem:unitriangular:basis}
Every family $\sY=\{y_v\mid v\in\B\}$ satisfying $y_v-v\in\V_{<v}$ is a
unitriangular basis. The transformation $T(v)=y_v$ and its inverse
belong to $\UTBV$, and
$\Span\{y_u\mid u\leq v\}=\V_{\leq v}$ for every $v\in\B$.
Consequently, $\UTBV=\UT_{\sY}(\V)$, and conjugation by $T$ on $\UTBV$ preserves
diagonal coordinates.
\end{lemma}
\begin{proof}
In any finite nonzero linear combination of the $y_v$, the coefficient
at the greatest index cannot cancel. Transfinite induction gives
$\Span\{y_u\mid u\leq v\}=\V_{\leq v}$.
Thus $T$ is invertible and both $T$ and $T^{-1}$ preserve these subspaces.
The equality of the triangular algebras follows. Diagonal coordinates
are preserved because $T(v,v)=T^{-1}(v,v)=1$ for every $v$.
\end{proof}

\begin{proposition}[{\cite[Proposition~6]{WY}}]\label{prop:corner:iso:ut}
Let $\ff\in\UTBV$ be an associative idempotent, let
$\B'=\Supp(\ff)$, put $\V'=\Span\B'$, and let $\pp$ be the diagonal
projection onto $\V'$. Then there is an invertible $T\in\UTBV$ with
$T(v,v)=1$ for every $v\in\B$ such that $\ff=T\pp T^{-1}$. Moreover,
$\range(\ff)=\Span\{\ff(v)\mid v\in\B'\}$, and
$\ff\UTBV\ff\cong\UT_{\B'}(\V')$ as topological $\F$-algebras.
The isomorphism preserves diagonal coordinates on $\B'$ and hence
preserves diagonal supports. In particular, a nonzero associative
idempotent has nonempty support.
\end{proposition}
\begin{proof}
Set $y_v=\ff(v)$ for $v\in\B'$ and $y_v=(1-\ff)(v)$ otherwise.
This family is unitriangular. By \Cref{lem:unitriangular:basis}, the
map $T(v)=y_v$ is invertible and triangular with $T(v,v)=1$, and
$\ff T=T\pp$, so $\ff=T\pp T^{-1}$. Hence
$\range(\ff)=T(\V')=\Span\{y_v\mid v\in\B'\}$, and $\ff=0$ if
$\B'=\emptyset$.
Conjugation by $T^{-1}$ identifies $\ff\UTBV\ff$ with
$\pp\UTBV\pp$, which restricts to $\UT_{\B'}(\V')$ on $\V'$.
Both maps and their inverses are continuous in the finite topology.
Conjugation preserves diagonal coordinates, and every element of
$\ff\UTBV\ff$ has zero diagonal outside $\B'$.
\end{proof}

\begin{lemma}\label{lem:finite:corners}
Let $\aA=\UTBVlim$ or $\aA=\UTBVfin$, and let $\ee\in\aA$ be an
associative idempotent. Then $\ee\aA\ee=\ee\UTBV\ee$.
\end{lemma}
\begin{proof}
For $\aA=\UTBVlim$, choose a finite subset $I\subseteq\B$ and
$\pp=\sum_{v\in I}\ee_{vv}$ such that $\ee=\pp\ee\pp$.
For every $x\in\UTBV$, we have
$\ee x\ee=\ee(\pp x\pp)\ee\in\ee\aA\ee$, since
$\pp x\pp\in\UTIV$.
For $\aA=\UTBVfin$, the inclusion follows from
$\range(\ee x\ee)\subseteq\range(\ee)$.
The reverse inclusion holds in both cases.
\end{proof}

\begin{lemma}\label{lem:idempotent:interval:support}
Let $\ee\in\UTBV$ be an associative idempotent. If $u\leq v$ and
$[u,v]\cap\Supp(\ee)=\emptyset$, then $\ee(u,v)=0$.
\end{lemma}
\begin{proof}
Write $\ee=T\pp T^{-1}$ as in \Cref{prop:corner:iso:ut}. Then
\[
\ee(u,v)=\sum_{w\in[u,v]\cap\Supp(\ee)}T(u,w)T^{-1}(w,v).
\]
Only finitely many summands can be nonzero, since $T^{-1}(v)$ has
finite support. The sum is zero under the stated hypothesis.
\end{proof}

\begin{lemma}\label{lem:singleton:idempotent:rankone}
Let $\ee\in\UTBV$ be an associative idempotent with
$\Supp(\ee)=\{v\}$, and set $y_\ee=\ee(v)$.
Then $\range(\ee)=\F y_\ee$, and $\ee(u)=0$ for $u<v$.
In particular, $\ee$ has rank one and is primitive in every subspace
of $\UTBV$ containing it.
\end{lemma}
\begin{proof}
The range formula is the one in \Cref{prop:corner:iso:ut}. For $u<v$,
the vector $\ee(u)\in\V_{<v}$ is a multiple of $y_\ee$, whose
coefficient at $v$ is $\ee(v,v)=1$. Hence $\ee(u)=0$. A rank-one
idempotent is primitive in every subspace containing it, since
nonzero orthogonal idempotents have independent ranges.
\end{proof}

\begin{lemma}\label{lem:half:orthogonal}
Let $\ee,\ff\in\UTBV$ be associative idempotents with
$\Supp(\ee)=\{u\}$. If every element of $\Supp(\ff)$ is greater
than $u$, then $\ff\ee=0$. If every element of $\Supp(\ff)$ is
less than $u$, then $\ee\ff=0$.
\end{lemma}
\begin{proof}
In the first case, \Cref{lem:idempotent:interval:support} shows that
$\ff$ vanishes on $\V_{\leq u}$, which contains $\range(\ee)$.
In the second case, the range formula in
\Cref{prop:corner:iso:ut} gives $\range(\ff)\subseteq\V_{<u}$,
and \Cref{lem:singleton:idempotent:rankone} shows that $\ee$ vanishes
on this subspace.
\end{proof}

\begin{lemma}\label{lem:unitriangular:basis:from:idempotents}
Suppose that $\ee_v\in\UTBV$ is an associative idempotent with
$\Supp(\ee_v)=\{v\}$ for every $v\in\B$. Then
$\sY=\{\ee_v(v)\mid v\in\B\}$ is a unitriangular basis of $\V$.
\end{lemma}
\begin{proof}
Idempotency gives $\ee_v(v,v)=1$, so
\Cref{lem:unitriangular:basis} applies.
\end{proof}

\begin{lemma}\label{lem:coordinate:projections:unitriangular:basis}
Let $\sY=\{y_v\mid v\in\B\}$ be a unitriangular basis, and let
$\qq_v$ be its coordinate projection onto $\F y_v$.
Then $\qq_v\in\UTBVfin$ is a primitive associative idempotent with
$\Supp(\qq_v)=\{v\}$. The family $\{\qq_v\mid v\in\B\}$ is
pairwise orthogonal and summable in the finite topology, with sum $1$.
\end{lemma}
\begin{proof}
By \Cref{lem:unitriangular:basis},
$\qq_v=T\ee_{vv}T^{-1}$ is triangular and has the same diagonal as
$\ee_{vv}$, where $T(v)=y_v$.
It has rank one, so it belongs to $\UTBVfin$ and is primitive.
Coordinate projections are pairwise orthogonal. Every finite-dimensional
subspace lies in the span of finitely many $y_v$, on which the
corresponding finite sums act as the identity. Thus their sum is $1$
in the finite topology.
\end{proof}

Let $\aA$ be an associative subalgebra of $\UTBV$, and put
$\range\aA=\sum_{x\in\aA}\range(x)$. A family $\sE\subseteq\aA$
of nonzero pairwise orthogonal associative idempotents is
\highlight{complete in $\aA$} if
\[
\range\aA=\bigoplus_{\ee\in\sE}\range(\ee).
\]
This holds precisely when
$\sum_{\ee\in\sE}\ee|_{\range\aA}=\id_{\range\aA}$ in the finite
topology. Indeed, orthogonality makes each $\ee$ the projection
onto its range, and convergence on each vector means that only
finitely many of these projections are nonzero there
(cf. \cite[Corollary~2.12]{IMesR}).
If $\aA$ contains every $\ee_{vv}$, then $\range\aA=\V$ and
completeness means $\sum_{\ee\in\sE}\ee=1$ in $\EndFV$.
For $\aA=\ff\UTBV\ff$, it means $\sum_{\ee\in\sE}\ee=\ff$,
since $\range\aA=\range(\ff)$ and every $\ee\in\aA$ vanishes on
$\ker(\ff)$.

\begin{definition}
Let $\aA$ be an associative subalgebra of \UTBV endowed with an $\Omega_\lambda$-grading. A
subset $\sE\subseteq\aA_1$ is called a \highlight{\sset} of
$\aA$ if it is complete in $\aA$ and consists of pairwise
orthogonal homogeneous associative idempotents whose supports are
singletons.
\end{definition}
For example, if $\aA=\UTBVlim$ is equipped with the elementary $\Omega_\lambda$-grading $\Gamma(\gamma)$, then $\deg(\ee_{vv})=\gamma(v)\gamma(v)^{-1}=1$ for all $v\in\B$ and $\{\ee_{vv}\mid v\in\B\}\subset\aA_1$ is a \sset of \aA.

\begin{lemma}\label{lem:sset:diagonalization}
Let $\aA\subseteq\UTBV$ be an associative subalgebra endowed with an $\Omega_\lambda$-grading, and let
$\sE=\{\ee_v\mid v\in\B\}\subseteq\aA_1$
be a \sset of \aA, indexed so that $\Supp(\ee_v)=\{v\}$.
Then there exists a well-ordered basis
$\B'=\{y_v\mid v\in\B\}$
of $\V$, ordered by
$y_u\leq y_v$ whenever $u\leq v$,
such that
$\UTBV=\UT_{\B'}(\V)$
and every element of $\sE$ is diagonal with respect to $\B'$.
More precisely,
$\ee_v=\ee_{y_vy_v}^{\,\B'}$
for every $v\in\B$, where $\ee_{y_vy_v}^{\,\B'}$ denotes the
corresponding diagonal matrix unit relative to $\B'$.
\end{lemma}

\begin{proof}
Set $y_v=\ee_v(v)$. By
\Cref{lem:unitriangular:basis:from:idempotents,lem:unitriangular:basis},
these vectors form a unitriangular basis and
$\UTBV=\UT_{\B'}(\V)$. Orthogonality gives
$\ee_v(y_w)=\ee_v\ee_w(w)=\delta_{v,w}y_w$.
Thus $\ee_v$ is the coordinate projection onto $\F y_v$.
\end{proof}

\begin{lemma}\label{lem:gen:class:sset}
Let $\UTBVlim\subseteq\aA\subseteq\UTBV$, and assume that $\aA$ is
an associative subalgebra endowed with an $\Omega_\lambda$-grading,
where $\lambda\notin\{0,-1\}$ and $\chr\F\neq2$ if $\lambda=1$.
If $\aA_1$ contains a \sset of
$\aA$, then the grading is elementary relative to a unitriangular
basis $\B'$ and is also an associative grading. If $T$ carries
$\B$ to $\B'$, conjugation by $T^{-1}$ carries the grading to an
elementary grading relative to $\B$ on $T^{-1}\aA T$.
For $\aA=\UTBVfin$ or $\aA=\UTBV$, we have
$T^{-1}\aA T=\aA$.
\end{lemma}

\begin{proof}
Let $\sE\subseteq\aA_1$ be a \sset of $\aA$. Since
$\range\aA=\V$, completeness gives
$\sum_{\ee\in\sE}\ee=1$. For each $v\in\B$, taking the
$(v,v)$-coordinate shows that some $\ee_v\in\sE$ has
$\Supp(\ee_v)=\{v\}$. It is unique, because diagonal coordinates
multiply, so two orthogonal idempotents cannot both have
$(v,v)$-coordinate $1$. Hence $\sE=\{\ee_v\mid v\in\B\}$.

Set $y_v=\ee_v(v)$. By \Cref{lem:sset:diagonalization},
$\B'=\{y_v\mid v\in\B\}$ is a well-ordered unitriangular basis
such that $\UTBV=\UT_{\B'}(\V)$ and $\ee_v$ is the projection
onto $\F y_v$. Let $\ee'_{uv}$ denote the matrix units
relative to $\B'$. Since
$\ee'_{uv}=\ee_u\ee_{uv}\ee_v$, we have
$\ee'_{uv}\in\aA$.

We claim that every $\ee'_{uv}$ is homogeneous. If $u=v$, then
$\ee'_{uu}=\ee_u\in\aA_1$. If $u<v$, write
$\ee'_{uv}=\sum_{g\in S}x_g$, where $x_g\in\aA_g$, and choose
$g$ such that $x_g(y_u,y_v)\ne0$. Then
$\ee'_{uu}\circ_\lambda
(\ee'_{vv}\circ_\lambda x_g)
=\lambda x_g(y_u,y_v)\ee'_{uv}\ne0$.
Thus $\ee'_{uv}$ is homogeneous of degree $g$.

Set $d(u,v)=\deg\ee'_{uv}$. For $u\leq v\leq w$, the element
$\ee'_{uv}\circ_\lambda\ee'_{vw}$ is a nonzero scalar
multiple of $\ee'_{uw}$. The scalar is $1+\lambda$ if
$u=v=w$, and $1$ otherwise. Hence
\begin{equation}\label{eq:gen:class:sset:cocycle}
d(u,v)d(v,w)=d(u,w).
\end{equation}
Let $v_0$ be the least element of $\B$ and define
$\gamma(y_v)=d(v_0,v)^{-1}$. By
\cref{eq:gen:class:sset:cocycle},
$d(u,v)=\gamma(y_u)\gamma(y_v)^{-1}$.

Let $x\in\aA_g$ be homogeneous. If $u<v$ and
$x(y_u,y_v)\ne0$, then
$\ee'_{uu}\circ_\lambda
(\ee'_{vv}\circ_\lambda x)
=\lambda x(y_u,y_v)\ee'_{uv}\ne0$, and hence $d(u,v)=g$.
If $u=v$ and $x(y_u,y_u)\ne0$, then
$\ee'_{uu}x\ee'_{uu}
=x(y_u,y_u)\ee'_{uu}\ne0$. This element is homogeneous of
degree $g$ by \Cref{product} if $\lambda\ne1$, and equals
$U_{\ee'_{uu}}(x)$ if $\lambda=1$. Thus $d(u,u)=g$ as well.
Therefore
$x(y_u,y_v)\ne0$ implies $\gamma(y_u)\gamma(y_v)^{-1}=g$.

Conversely, suppose that $x\in\aA$ satisfies this implication and
write $x=\sum_{h\in S}x_h$, where $x_h\in\aA_h$. A nonzero
coordinate of $x_h$ can occur only at a pair $(u,v)$ for which
$d(u,v)=h$. Thus the coordinate supports of the $x_h$ are
disjoint, and the hypothesis on $x$ implies that $x_h=0$ for
every $h\ne g$. Consequently,
\[
\aA_g=
\left\{
x\in\aA\ \middle|\
x(y_u,y_v)\ne0
\Longrightarrow
\gamma(y_u)\gamma(y_v)^{-1}=g
\right\}.
\]
Hence the grading is elementary relative to $\B'$.

This grading is associative. Indeed, for $x\in\aA_g$ and
$z\in\aA_h$, each nonzero summand in
$(xz)(y_u,y_v)=\sum_w x(y_u,y_w)z(y_w,y_v)$
satisfies $d(u,w)d(w,v)=gh$. The sum is finite because $z(y_v)$
is a finite linear combination. Thus $xz\in\aA_{gh}$.

Finally, let $T(v)=y_v$ and set
$\bar\gamma(v)=\gamma(y_v)$. Then
$T^{-1}\ee'_{uv}T=\ee_{uv}$, so conjugation by $T^{-1}$
carries the grading onto the elementary grading
$\Gamma(\bar\gamma)$ on $T^{-1}\aA T$.
Conjugation preserves $\UTBV$ by \Cref{lem:unitriangular:basis}
and preserves $\UTBVfin$ because it preserves rank.
\end{proof}

The preceding lemma reduces the classification problem to elementary
gradings. We now show that, for algebras containing $\UTBVfin$, every
grading elementary relative to $\B$ has finite support.

\begin{lemma}[{\cite[Lemma~22]{WY}}]\label{lem:finite:grading}
Let $\Gamma$ be a $G$-grading on $\aA$, where
$\UTBVfin\subseteq\aA\subseteq\UTBV$, which is elementary relative to
$\B$ with defining map $\gamma$. Then $\Supp\Gamma$ is finite.
\end{lemma}

\begin{proof}
Let $v_0$ be the least element of $\B$, and define
$x\in\UTBVfin$ by $x(v)=v_0$ for every $v\in\B$.
Write $x=\sum_{g\in S}x_g$, where $x_g\in\aA_g$ and $S$ is finite.
Since $x(v_0,v)=1$ and $\Gamma$ is elementary relative to $\B$,
$\gamma(v_0)\gamma(v)^{-1}\in S$ for every $v\in\B$.
Hence $\gamma(\B)\subseteq S^{-1}\gamma(v_0)$ is finite, and so is
$\Supp\Gamma\subseteq\gamma(\B)\gamma(\B)^{-1}$.
\end{proof}

\section{Gradings on finite triangular algebras}

Throughout this section, $\lambda\in\F\setminus\{0,\pm1\}$.

\begin{corollary}\label{grad_is_grad}
Let $G$ be abelian. A vector space $G$-grading on an associative
algebra is an $\Omega_\lambda$-grading if and only if it is an
associative grading. Two such gradings are isomorphic as
$\Omega_\lambda$-algebras if and only if they are isomorphic as
associative algebras.
\end{corollary}
\begin{proof}
The first assertion follows from \Cref{product} and the commutativity
of $G$. The second follows from \Cref{hom_is_hom}.
\end{proof}
Thus, for abelian $G$, the two classification problems coincide.

Fix a nonempty finite subset $I\subseteq\B$. We now turn to
$\Omega_\lambda$-gradings on $\UTIV$ by arbitrary groups. We begin with
a few technical lemmas.

\begin{lemma}\label{keylemma}
Let $\Gamma$ be an $\Omega_\lambda$-grading on $\UTIV$. Let $u$ and $v$
be the least and greatest elements of $I$, respectively. Then there
exist $x_u,x_v\in J(\UTIV)$ such that $\ee_{uu}+x_u$ and
$\ee_{vv}+x_v$ are homogeneous of degree~$1$.
\end{lemma}
\begin{proof}
Set $\aC=\UTIV$ and $\iJ=J(\aC)$. If $|I|=1$, then
$\ee_{uu}=1\in\aC_1$ by \Cref{lem:identity:neutral}, so assume that
$|I|>1$.

By \Cref{Jgraded}, $\iJ$ is graded. Define $\mathcal Z_1=\iJ$ and, recursively,
\[
\mathcal Z_{k+1}=\Span\{a\circ_\lambda b\mid a\in\mathcal Z_k,
b\in\iJ\}.
\]
Each $\mathcal Z_k$ is graded. A computation with matrix units gives
$\mathcal Z_k=\iJ^k$ for $1\leq k<|I|$, where the right-hand side is the
associative ideal power: inductively, every product lies in
$\iJ^{k+1}$, while
$\ee_{ab}\circ_\lambda\ee_{bc}=\ee_{ac}$ for $a<b<c$ supplies every
matrix unit in $\iJ^{k+1}$. Hence
$\mathcal Z_{|I|-1}=\Span\{\ee_{uv}\}$. Put $z=\ee_{uv}$. Then $z$ is
homogeneous, say of degree $d$.

For every homogeneous $a\in\aC_g$, triangularity gives
\[
a\circ_\lambda z=(a(u,u)+\lambda a(v,v))z
\quad\text{and}\quad
z\circ_\lambda a=(a(v,v)+\lambda a(u,u))z.
\]
If $g\ne1$, both products vanish by comparison of degrees.
Since $\lambda^2\ne1$, we obtain $a(u,u)=a(v,v)=0$.

Identify $\bar\aC=\aC/\iJ$ with $\F^I$ by the diagonal map.
Its associative multiplication is $(1+\lambda)^{-1}\circ_\lambda$
and is graded. Write $\bar\ee=\ee_{uu}+\iJ=\sum_h b_h$ with
$b_h\in\bar\aC_h$. The preceding observation gives
$b_1(u)=1$ and $b_g(u)=0$ for $g\ne1$.
Consequently, $b_1\bar\ee=\bar\ee$ and $b_g\bar\ee=0$.
Comparison of their degree-$g$ components gives
$b_1b_g=b_g$ and $b_gb_1=0$, respectively.
The quotient is commutative, so $b_g=0$ for $g\ne1$.
Thus $\ee_{uu}+\iJ$ is homogeneous of degree $1$.
The same argument applies to $\ee_{vv}+\iJ$.
Since $\iJ$ is graded, both cosets have representatives in $\aC_1$.
\end{proof}

\begin{lemma}\label{consequence_keylemma}
Let $\Gamma$ be an $\Omega_\lambda$-grading on $\UTIV$, and let $u$ and
$v$ be the least and greatest elements of $I$, respectively. Then there
is an invertible $S\in\UTIV$ such that $\ee_{uu}$ and $\ee_{vv}$ are
homogeneous of degree $1$ in the grading transported by
$x\mapsto S^{-1}xS$.
\end{lemma}
\begin{proof}
If $|I|=1$, then $\ee_{uu}=\ee_{vv}=1$ is homogeneous of degree $1$
by \Cref{lem:identity:neutral}, and $S=1$ works. Hence assume $|I|>1$.

By \Cref{keylemma}, there are $a_u,a_v\in(\UTIV)_1$ with
$a_u-\ee_{uu},a_v-\ee_{vv}\in J(\UTIV)$.
Put $n=|I|$. Each has characteristic polynomial $t^{n-1}(t-1)$
as a transformation of $\V_I=\Span I$.
By the Cayley-Hamilton theorem, $a_u^n=a_u^{n-1}$ and
$a_v^n=a_v^{n-1}$. Thus
$\ee_u=a_u^{n-1}$ and $\ee_v=a_v^{n-1}$ are associative idempotents
with supports $\{u\}$ and $\{v\}$, respectively.
They belong to $(\UTIV)_1$ by \Cref{product}.
By \Cref{lem:half:orthogonal}, $\ee_v\ee_u=0$.
Replacing $\ee_u$ by $\ee_u-\ee_u\ee_v$ as in
\Cref{lem:orthogonalize}, we obtain two orthogonal homogeneous
associative idempotents with the same supports.

For $t\in I\setminus\{u,v\}$, set
$y_t=(1-\ee_u-\ee_v)(t)$, and set
$y_u=\ee_u(u)$ and $y_v=\ee_v(v)$. The family
$\{y_t\mid t\in I\}$ is unitriangular relative to $I$, hence is a basis
of $\V_I$. Orthogonality gives
\[
\ee_u(y_t)=\delta_{u,t}y_u
\quad\text{and}\quad
\ee_v(y_t)=\delta_{v,t}y_v.
\]
Let $S(t)=y_t$ for $t\in I$. Then $S\in\UTIV$ is invertible, and these
identities give $S^{-1}\ee_uS=\ee_{uu}$ and $S^{-1}\ee_vS=\ee_{vv}$.
Since $\ee_u,\ee_v\in(\UTIV)_1$, both matrix units are homogeneous of
degree $1$ in the grading transported by $x\mapsto S^{-1}xS$.
\end{proof}

We are now in a position to classify the $\Omega_\lambda$-gradings on
$\UTIV$.
\begin{theorem}\label{nonLie_nonJordan_fd}
Let $\lambda\in\F\setminus\{0,\pm1\}$, and let $\Gamma$ be a group
$\Omega_\lambda$-grading on $\UTIV$. Then there exists an invertible
element $T\in\UTIV$ such that the grading
transported by the inner automorphism
$x\mapsto T^{-1}xT$ is elementary.
\end{theorem}

\begin{proof}
By \Cref{consequence_keylemma}, we may assume that $\ee_{uu}$ and
$\ee_{vv}$ are homogeneous of degree $1$, where $u$ and $v$ are the
least and greatest elements of $I$, respectively.

If $|I|\leq2$, then all the diagonal matrix units are homogeneous,
and the result follows from \Cref{lem:gen:class:sset}, applied with
ambient vector space $\V_I=\Span I$ and ordered basis $I$. Assume that
$|I|>2$, and set $I'=I\setminus\{u,v\}$. The subspace
\[
\aB=(1-\ee_{uu}-\ee_{vv})\UTIV(1-\ee_{uu}-\ee_{vv})
\]
is graded: indeed, $1\in(\UTIV)_1$ by
\Cref{lem:identity:neutral}, so this follows from
\Cref{rem:grad:inheritance}. Moreover,
$\aB=\UT_{I'}(\Span I')$. By induction, there
exists an invertible triangular transformation $S'$ of $\Span I'$
whose conjugation makes the grading on $\aB$ elementary.

Extend $S'$ to $\Span I$ by setting $S(u)=u$, $S(v)=v$, and
$S|_{\Span I'}=S'$. Then $S\in\UTIV^\times$, and the grading
transported by conjugation by $S$ has all the diagonal matrix units
homogeneous. Since these form a \sset in $(\UTIV)_1$, the result follows
from \Cref{lem:gen:class:sset}, again with ambient vector space
$\V_I$ and ordered basis $I$.
\end{proof}

\section{Gradings on infinite triangular algebras}

For the finite-map and finite-range algebras, homogeneous idempotents
can be constructed from algebraic elements.

\begin{lemma}\label{existencehomogeneous}
Let $\aA=\UTBVlim$ or $\aA=\UTBVfin$, endowed with an
$\Omega_\lambda$-grading, where either
$\lambda\in\F\setminus\{0,\pm1\}$, or $\lambda=1$ and
$\chr\F\neq2$. For every $v\in\B$, there exists an associative
idempotent $\ff\in\aA_1$, primitive in $\aA_1$, such that
$v\in\Supp(\ff)$.
\end{lemma}

\begin{proof}
Choose a homogeneous component $x$ of $\ee_{vv}$ with $x(v,v)\neq0$.
Since $x$ has finite range, there is a finite $I\subseteq\B$ such
that $\range(x)\subseteq\Span I$. The restriction of $x$ to
$\Span I$ is upper triangular. If $q_0$ is its characteristic
polynomial, then $q_0(x)x=0$ on $\V$. Thus $x$ is algebraic.

For every $m\ge1$, we have
\[
x^{\circ_\lambda m}=(1+\lambda)^{m-1}x^m\neq0.
\]
These powers span a finite-dimensional subspace, so $\deg x$ has
finite order. Replacing $x$ by the corresponding associative power,
we may assume that $x\in\aA_1$. All its positive associative powers
then belong to $\aA_1$.

Write the minimal polynomial as $m_x(t)=t^r h(t)$, where $h(0)\neq0$.
The Chinese remainder theorem yields $p(t)\in\F[t]$ such that
\[
p(t)\equiv0\pmod{t^{\max\{1,r\}}}
\quad\text{and}\quad
p(t)\equiv1\pmod{h(t)}.
\]
Then $\ee=p(x)\in\aA_1$ is an associative idempotent, since
$m_x$ divides $p^2-p$ and $p(0)=0$. Since $x(v,v)\neq0$ is a root of $h$,
$\ee(v,v)=p(x(v,v))=1$. If $\ee$ is not primitive in $\aA_1$,
decompose it into two nonzero orthogonal idempotents in $\aA_1$
and retain the summand whose support contains $v$. The summands
have disjoint supports, which are nonempty by
\Cref{prop:corner:iso:ut}. Since $\Supp(\ee)$ is finite, this process
terminates.
\end{proof}

\begin{lemma}\label{lem:primitive:support:finite:corner}
Let $\aA=\UTBVlim$, $\UTBVfin$, or $\UTBV$, endowed with an
$\Omega_\lambda$-grading, where either
$\lambda\in\F\setminus\{0,\pm1\}$, or $\lambda=1$ and
$\chr\F\neq2$. Let $\ee\in\aA_1$ be an associative idempotent
that is primitive in $\aA_1$ and has finite support. Then
$|\Supp(\ee)|=1$ if $\lambda\neq1$, and
$|\Supp(\ee)|\in\{1,2\}$ if $\lambda=1$.
\end{lemma}

\begin{proof}
Put $I=\Supp(\ee)$ and $n=|I|$. By
\Cref{rem:grad:inheritance,prop:corner:iso:ut},
the corner $\ee\aA\ee$ inherits the grading. It equals
$\ee\UTBV\ee$, trivially for $\aA=\UTBV$ and by
\Cref{lem:finite:corners} otherwise, and is isomorphic to
$\UT_I(\Span I)\cong\UT_n(\F)$. Its identity is
primitive in its neutral component.

For $\lambda\notin\{0,\pm1\}$, \Cref{nonLie_nonJordan_fd} makes
the transported grading elementary. If $n>1$, a diagonal matrix
unit and its complement contradict this primitivity. Thus $n=1$.
For $\lambda=1$, \cite[Theorem~24]{KY2017} gives an elementary or
an MT grading.\footnote{Although the introduction of \cite{KY2017}
mentions infinite fields, the normalization in Lemmas~19--23 uses
only $\chr\F\neq2$.}
If $n>2$, a diagonal matrix unit in the elementary case, or the sum
of the first and last diagonal matrix units in the MT case, is a
nontrivial neutral idempotent. Its complement again contradicts
primitivity. Hence $n\leq2$.
\end{proof}

\begin{corollary}\label{primitivehomogeneouselement}
Let $\lambda\in\F\setminus\{0,\pm1\}$, and let
$\aA=\UTBVlim$ or $\aA=\UTBVfin$ be endowed with an
$\Omega_\lambda$-grading. For every $v\in\B$, there exists an
associative idempotent $\ee_v\in\aA_1$ with $\Supp(\ee_v)=\{v\}$,
primitive in $\aA_1$ and in $\aA$.
\end{corollary}

\begin{proof}
Apply \Cref{existencehomogeneous,lem:primitive:support:finite:corner}.
Primitivity in $\aA$ follows from
\Cref{lem:singleton:idempotent:rankone}.
\end{proof}

\begin{lemma}\label{lem:primitive:support:orthogonal}
Assume that $\chr\F\neq2$, and let
$\aA=\UTBVlim$, $\UTBVfin$, or $\UTBV$ be endowed with an
$\Omega_1$-grading. Let $(\ee_i)_{i\geq1}$ be an infinite sequence
of pairwise orthogonal associative idempotents in $\aA_1$, each
primitive in $\aA_1$ and having finite support. Then every
$\Supp(\ee_i)$ is a singleton.
\end{lemma}

\begin{proof}
By \Cref{lem:primitive:support:finite:corner}, all supports have
one or two elements, and they are pairwise disjoint. After
reindexing, suppose that $|\Supp(\ee_1)|=2$.
For any finite sum $\ff$ of the $\ee_i$, the graded corner
$\ff\aA\ff$ is support-preservingly isomorphic to
$\UT_{\Supp(\ff)}(\Span\Supp(\ff))$ by
\Cref{rem:grad:inheritance,prop:corner:iso:ut}, using
\Cref{lem:finite:corners} when $\aA\neq\UTBV$.
Each summand remains primitive in the neutral component of this
corner, since any decomposition there is a decomposition in
$\aA_1$.

We use the following consequences of
\cite[Lemma~14 and the proofs of Lemmas~20--23]{KY2017} for a
grading on $\UT_4(\F)^{(+)}$, indexed by $s_1<s_2<s_3<s_4$.
Conjugation by an invertible triangular matrix carries it to an
elementary or an MT grading and preserves diagonal supports.
In the elementary case, the neutral component is the direct sum
of upper triangular algebras on the fibres of the defining map.
Its primitive associative idempotents therefore have singleton
supports. An MT grading has a homogeneous basis each of whose
elements is strictly triangular or equal to one of
$\ee_{s_1s_1}\pm\ee_{s_4s_4}$ and $\ee_{s_2s_2}\pm\ee_{s_3s_3}$, so
the diagonal part of a neutral element is neutral. The two
differences are homogeneous of a common degree of order $2$, hence
not neutral. Thus the diagonal part of the neutral component is
\[
\Span\{\ee_{s_1s_1}+\ee_{s_4s_4},
        \ee_{s_2s_2}+\ee_{s_3s_3}\}.
\]
Consequently, supports of nonzero neutral idempotents are unions
of the two mirror pairs.

There is at most one singleton support among the $\Supp(\ee_i)$.
Otherwise, the corner determined by $\ee_1$ and two idempotents
with singleton supports is isomorphic to $\UT_4(\F)^{(+)}$, and
its identity is a sum of primitive neutral idempotents with
support sizes $2,1,1$. The elementary case excludes the first
summand, and the MT case excludes the other two.

Infinitely many supports therefore have size two. Write them as
$\{u_i<v_i\}$. The corner determined by any two corresponding
idempotents is isomorphic to $\UT_4(\F)^{(+)}$. It cannot have an
elementary grading, so their supports are its mirror pairs.
Thus $u_i<u_j$ implies $v_i>v_j$. The infinite set of lower
endpoints in the well-ordered set $\B$ contains a strictly
increasing sequence. The corresponding upper endpoints form a
strictly decreasing sequence, a contradiction.
\end{proof}

\begin{lemma}\label{lem:primitive:support:infinite:jordan}
Assume that $\chr\F\neq2$ and that $\B$ is infinite. Let
$\aA=\UTBVlim$ or $\aA=\UTBVfin$, endowed with an
$\Omega_1$-grading. Every associative idempotent in $\aA_1$ that
is primitive in $\aA_1$ has singleton support.
\end{lemma}

\begin{proof}
It suffices, by \Cref{lem:primitive:support:orthogonal}, to
include the given idempotent in an infinite sequence of pairwise
orthogonal idempotents primitive in $\aA_1$.
Suppose that $\ee_1,\ldots,\ee_n$ have been chosen and put
$\ff=\ee_1+\cdots+\ee_n$. Choose
$v\in\B\setminus\Supp(\ff)$ and a homogeneous component $y$ of
$\ee_{vv}$ such that $y(v,v)\neq0$. By
\Cref{rem:grad:inheritance}, the element
\[
z=y-\ff\circ y+U_\ff(y)=(1-\ff)y(1-\ff)
\]
is homogeneous of the same degree as $y$. The last expression
is computed in $\EndFV$. It gives $\ff z=z\ff=0$ and
$z(v,v)=y(v,v)\neq0$.

The polynomial construction in the proof of
\Cref{existencehomogeneous}, applied to $z$, gives a neutral
idempotent $\pp$ with $v\in\Supp(\pp)$ and
$\ff\pp=\pp\ff=0$. Choose a primitive neutral summand
$\ee_{n+1}$ of $\pp$ whose support contains $v$. Since
$\ee_{n+1}=\pp\ee_{n+1}=\ee_{n+1}\pp$, it remains orthogonal to
$\ff$. This continues the sequence.
\end{proof}

\begin{lemma}\label{lem:finite:diagonalization}
Let $\UTBVfin\subseteq\aA\subseteq\UTBV$ be an associative
algebra endowed with an $\Omega_\lambda$-grading, where either
$\lambda\in\F\setminus\{0,\pm1\}$, or $\lambda=1$ and
$\chr\F\neq2$. Suppose that, for each $v\in\B$, there is an
associative idempotent $\ee_v\in\aA_1$ with
$\Supp(\ee_v)=\{v\}$. Set $y_v=\ee_v(v)$, and let $\qq_v$
be the coordinate projection onto $\F y_v$ relative to the basis
$\{y_w\mid w\in\B\}$. Then $\qq_v\in\UTBVfin\cap\aA_1$,
$\Supp(\qq_v)=\{v\}$, and $\{\qq_v\mid v\in\B\}$ is a
\sset\ of $\aA$.
\end{lemma}

\begin{proof}
By \Cref{lem:unitriangular:basis:from:idempotents,lem:coordinate:projections:unitriangular:basis}, the $y_v$ form
a unitriangular basis, and the $\qq_v$ are pairwise orthogonal
rank-one idempotents in $\UTBVfin$ with singleton supports and
sum $1$. It remains to prove that they belong to $\aA_1$.

For each nonempty finite $F=\{v_1<\cdots<v_n\}\subseteq\B$, construct
idempotents by descending induction. Put $\qq_{v_n}^F=\ee_{v_n}$
and, for $k<n$, set
\[
\pp_k=\sum_{j>k}\qq_{v_j}^F
\quad\text{and}\quad
\qq_{v_k}^F=\ee_{v_k}-\ee_{v_k}\pp_k
              =\ee_{v_k}-\ee_{v_k}\circ_\lambda\pp_k.
\]
At each step $\pp_k\ee_{v_k}=0$ by
\Cref{lem:half:orthogonal}. Hence \Cref{lem:orthogonalize}
shows that the $\qq_{v_k}^F$ are pairwise orthogonal idempotents
in $\aA_1$ with unchanged singleton supports and
$\range(\qq_{v_k}^F)=\F y_{v_k}$.
The same induction gives
\[
\qq_u^F(y_w)=\delta_{uw}y_u\qquad(u,w\in F).
\]
Indeed, $\pp_k$ fixes $y_{v_j}$ for $j>k$ and annihilates
$\V_{\leq v_k}$, while $\ee_{v_k}$ fixes $y_{v_k}$ and
annihilates $\V_{<v_k}$.

Fix $v\in\B$ and write $\qq_v=\sum_g z_g$ homogeneously.
For every finite $F$ containing $v$, put
$\pp_F=\sum_{u\in F}\qq_u^F$. Its range is
$\Span\{y_u\mid u\in F\}$, and the preceding identities give
\[
\pp_F\qq_v\pp_F=\qq_v^F\in\aA_1.
\]
By \Cref{rem:grad:inheritance}, each $\pp_Fz_g\pp_F$ belongs
to $\aA_g$, so $\pp_Fz_g\pp_F=0$ whenever $g\neq1$.
For any $w\in\B$, choose $F$ containing $v,w$ and the finite
support of $z_g(y_w)$ in the $y$-basis. Then $\pp_F$ fixes both
$y_w$ and $z_g(y_w)$, whence $z_g(y_w)=0$. Thus $z_g=0$ for
$g\neq1$, as required. Since $\range\aA=\V$, the coordinate
projections form a \sset\ of $\aA$.
\end{proof}

\begin{lemma}\label{lem:existence:sset:intermediate}
Let $\aA$ be an associative subalgebra with
$\UTBVfin\subseteq\aA\subseteq\UTBV$, endowed with an
$\Omega_\lambda$-grading, where either
$\lambda\in\F\setminus\{0,\pm1\}$, or $\lambda=1$, $\chr\F\neq2$,
and $\B$ is infinite. If $\UTBVfin$ is a graded ideal of $\aA$,
then $\aA_1$ contains a \sset
$\{\qq_v\mid v\in\B\}$ with
$\qq_v\in\UTBVfin\cap\aA_1$ and $\Supp(\qq_v)=\{v\}$.
\end{lemma}

\begin{proof}
Apply \Cref{primitivehomogeneouselement} to the induced grading
on $\UTBVfin$ if $\lambda\neq1$, and
\Cref{existencehomogeneous,lem:primitive:support:infinite:jordan}
if $\lambda=1$. These results give neutral idempotents with
prescribed singleton supports. \Cref{lem:finite:diagonalization}
then gives the required \sset.
\end{proof}

\begin{lemma}\label{lem:countable:half:orthogonal:straightening}
Let $\B=\{v_n\mid n\in\NN\}$ with $v_1<v_2<\cdots$, and let
$\aA=\UTBVlim$ be endowed with an $\Omega_\lambda$-grading, where
either $\lambda\in\F\setminus\{0,\pm1\}$, or $\lambda=1$ and
$\chr\F\neq2$. Suppose that for each $n$ there is an
associative idempotent $\ee_n\in\aA_1$ with
$\Supp(\ee_n)=\{v_n\}$. Then $\aA_1$ contains
a \sset\ $\{\ff_n\mid n\in\NN\}$ with
$\Supp(\ff_n)=\{v_n\}$.
\end{lemma}

\begin{proof}
Set $\ff_1=\ee_1$. Having constructed pairwise orthogonal
idempotents $\ff_1$, \dots, $\ff_{n-1}\in\aA_1$ with the required
supports, put $\pp_{n-1}=\sum_{i<n}\ff_i$ and
\[
\ff_n=\ee_n-\pp_{n-1}\ee_n
      =\ee_n-\pp_{n-1}\circ_\lambda\ee_n.
\]
By \Cref{lem:half:orthogonal}, $\ee_n\pp_{n-1}=0$.
Thus \Cref{lem:orthogonalize} gives a neutral idempotent
$\ff_n$ orthogonal to $\pp_{n-1}$ with
$\Supp(\ff_n)=\{v_n\}$. It is primitive by
\Cref{lem:singleton:idempotent:rankone}.

By \Cref{lem:unitriangular:basis:from:idempotents}, the vectors
$\ff_n(v_n)$ form a basis of $\V$. Their spans are the ranges
of the pairwise orthogonal $\ff_n$, so the family is complete.
\end{proof}

\begin{lemma}\label{lem:countable:finite:maps}
Suppose that $\B$ is countably infinite. Every associative group
grading on $\UTBVlim$ is carried to an elementary grading relative
to $\B$ by conjugation by an invertible triangular transformation
that normalizes $\UTBVlim$.
\end{lemma}

\begin{proof}
Put $\aA=\UTBVlim$, and enumerate its matrix units as
$x_1$, $x_2$, \dots{} We construct neutral associative idempotents
$\ff_n$, starting with $\ff_0=0$, such that
\[
\ff_n\ff_{n-1}=\ff_{n-1}=\ff_{n-1}\ff_n
\quad\text{and}\quad
\ff_nx_i=x_i=x_i\ff_n\qquad(i\leq n).
\]
Let $X_n$ consist of $\ff_{n-1}$ and all homogeneous components
of $x_1$, \dots, $x_n$. Choose a finite diagonal idempotent
$\pp\in\aA$ satisfying $\pp x=x=x\pp$ for $x\in X_n$.
If $a$ is the neutral component of $\pp$, comparison of degrees
gives $ax=x=xa$ for $x\in X_n$. The finite map $a$ is algebraic.
Write $m_a(t)=t^rh(t)$, where $h(0)\neq0$. Since $X_n$ contains
a nonzero element, $m_a(1)=0$, so $h(1)=0$. The Chinese remainder
theorem gives $q(t)$ satisfying
\[
q(t)\equiv0\pmod{t^{\max\{1,r\}}}
\quad\text{and}\quad
q(t)\equiv1\pmod{h(t)}.
\]
Then $\ff_n=q(a)\in\aA_1$ is an associative idempotent.
Since $q(0)=0$ and $q(1)=1$, it satisfies the required identities.

Each nonzero $\ff_n-\ff_{n-1}$ is a neutral associative
idempotent. By \Cref{lem:finite:corners,prop:corner:iso:ut} and
the finite associative classification in \cite{VZ2007}, it is
a sum of orthogonal neutral associative idempotents with
singleton supports. Indeed, each primitive summand has a
one-dimensional corner, whereas \Cref{prop:corner:iso:ut}
identifies that corner with the triangular algebra on its support.
These decompositions give pairwise orthogonal
idempotents $\{\ee_v\mid v\in\B\}\subseteq\aA_1$ with
$\Supp(\ee_v)=\{v\}$, since every $\ee_{vv}$ has a local unit
among the $\ff_n$.

Set $T(v)=\ee_v(v)$. By
\Cref{lem:unitriangular:basis:from:idempotents}, $T$ is invertible
and triangular. Every $x\in\aA$ satisfies $x=\ff_nx\ff_n$ for
some $n$, so it is a finite sum of terms $\ee_u x\ee_v$.
Every nonzero term has $u\leq v$ and is a scalar multiple of
$T\ee_{uv}T^{-1}$.
Conversely, $T\ee_{uv}T^{-1}=\ee_u\ee_{uv}\ee_v\in\aA$.
Hence
\[
\aA=\Span\{T\ee_{uv}T^{-1}\mid u\leq v\},
\]
so $T^{-1}\aA T=\aA$. The transported grading has neutral
diagonal matrix units. Its one-dimensional corners
$\F\ee_{uv}$ are therefore graded. Their degrees satisfy
$d(u,v)d(v,w)=d(u,w)$, so taking
$\gamma(v)=d(v_0,v)^{-1}$, where $v_0=\min\B$, proves that the
grading is elementary.
\end{proof}

\begin{theorem}\label{thm:classification:unified}
Let $\Gamma$ be an $\Omega_\lambda$-grading on an associative
subalgebra $\aA\subseteq\UTBV$, where either
$\lambda\in\F\setminus\{0,\pm1\}$, or $\lambda=1$, $\chr\F\neq2$,
and $\B$ is infinite. Assume that one of the following conditions
holds:
\begin{enumerate}
\item\label{thm:classification:unified:1}
$\UTBVfin\subseteq\aA\subseteq\UTBV$ and $\UTBVfin$ is a graded
ideal of $\aA$.
\item\label{thm:classification:unified:2}
$\aA=\UTBVlim$ and $\B=\{v_n\mid n\in\NN\}$, where
$v_1<v_2<\cdots$.
\end{enumerate}
Then $\Gamma$ is elementary relative to a unitriangular basis of
$\V$ and is also an associative grading. If $\aA=\UTBVfin$,
$\aA=\UTBV$, or condition \labelcref{thm:classification:unified:2}
holds, then $\Gamma$ is isomorphic to an elementary grading on
$\aA$ relative to $\B$.
\end{theorem}

\begin{proof}
Under \labelcref{thm:classification:unified:1},
\Cref{lem:existence:sset:intermediate} gives a \sset\ in $\aA_1$.
Under \labelcref{thm:classification:unified:2},
\Cref{primitivehomogeneouselement} if $\lambda\neq1$, and
\Cref{existencehomogeneous,lem:primitive:support:infinite:jordan}
if $\lambda=1$, give neutral idempotents with prescribed singleton
supports. Apply \Cref{lem:countable:half:orthogonal:straightening}.
In either case, \Cref{lem:gen:class:sset} proves the first assertions.
For the last assertion, that lemma applies directly to
$\UTBVfin$ and $\UTBV$. Under
\labelcref{thm:classification:unified:2}, apply
\Cref{lem:countable:finite:maps} to the associative grading just
obtained. Associative isomorphisms preserve $\circ_\lambda$.
\end{proof}

\begin{corollary}
Let $\Gamma$ be an $\Omega_1$-grading on $\UTBVlim$, where
$\chr\F\neq2$ and $\B=\{v_n\mid n\in\NN\}$ with
$v_1<v_2<\cdots$. Then $\Supp\Gamma$ generates an abelian group.
\end{corollary}

\begin{proof}
By \Cref{thm:classification:unified}, the grading is isomorphic
to an elementary grading on the same algebra.
Apply \Cref{prop:elementary-omega}.
\end{proof}

\begin{corollary}
Let $\Gamma$ be an $\Omega_\lambda$-grading on an associative
subalgebra $\aA$ such that
\[
\UTBVfin\subseteq\aA\subseteq\UTBV,
\]
and $\UTBVfin$ is a graded ideal of $\aA$. Assume that either
$\lambda\in\F\setminus\{0,\pm1\}$, or $\lambda=1$ and
$\chr\F\ne2$. Then $\Supp\Gamma$ is finite.
\end{corollary}

\begin{proof}
If $\B$ is finite, then $\aA$ is finite-dimensional, so the
support is finite. Otherwise, \Cref{thm:classification:unified}
makes the grading elementary relative to a unitriangular basis.
The conjugation in \Cref{lem:gen:class:sset} preserves $\UTBVfin$, so the
transported algebra still contains $\UTBVfin$. Its support is
finite by \Cref{lem:finite:grading}.
\end{proof}

\section{Closed gradings on the full triangular algebra}
We now consider closed $\Omega_\lambda$-gradings on $\UTBV$, without
assuming that $\UTBVfin$ is graded. We construct homogeneous
idempotents by passing to the diagonal quotient and lifting them.

\begin{proposition}\label{lem:jacobson:radical}
Let $\lambda\in\F^\times$, and endow $\UTBV$ with an
$\Omega_\lambda$-grading. If $\aA\subseteq\UTBV$ is a graded
associative subalgebra satisfying
$J(\aA)=\aA\cap J(\UTBV)$, then $J(\aA)$ is graded.
\end{proposition}

\begin{proof}
It suffices to prove the assertion for $\aA=\UTBV$.
Write $[a,b]=ab-ba$, and let $\sS$ be the span of all
$\Omega_\lambda$-associators in $\aA$. The identity
\[
(a,b,c)_\lambda=\lambda[b,[a,c]]
\]
shows that every associator is strictly triangular. Hence
$\sS\subseteq J(\aA)$ by \Cref{lem:radicals:finite-algebras}.

Let $x\in J(\aA)$. Define $S\in\aA$ by sending each basis vector
to its immediate predecessor, or to zero if it has none.
Let $R\in\EndFV$ send each basis vector to its immediate successor,
or to zero if it has none. For every $v\in\B$, we have
$R(\V_{<v})\subseteq\V_{\leq v}$ and
$SR|_{\V_{<v}}=\id_{\V_{<v}}$.
Recursively over $\B$, define
\[
Y(v)=R\bigl(x(v)+Y(S(v))\bigr),
\]
where $Y(0)=0$. Inductively, the vector inside $R$ belongs to
$\V_{<v}$.
Thus this defines a triangular transformation $Y\in\aA$ and gives
$SY(v)=x(v)+YS(v)$ for every $v$, so $x=[S,Y]$.

Choose $\varepsilon(v)\in\{0,1\}$ recursively, taking
$\varepsilon(v)=0$ when $v$ has no immediate predecessor and
$\varepsilon(v)=1-\varepsilon(S(v))$ otherwise.
Let $\ee$ be the diagonal associative idempotent defined by
$\ee(v)=\varepsilon(v)v$, and put $T=[\ee,S]$.
Whenever $S(v)\neq0$, we have
$(\varepsilon(S(v))-\varepsilon(v))^2=1$.
Consequently $[\ee,T]=S$ in every characteristic, and
\[
\begin{aligned}
x=[S,Y]
&=[\ee,[T,Y]]-[T,[\ee,Y]]\\
&=\lambda^{-1}\bigl((T,\ee,Y)_\lambda-(\ee,T,Y)_\lambda\bigr)
\in\sS.
\end{aligned}
\]
Therefore $J(\aA)=\sS$. If $a$, $b$, $c$ have degrees $g$, $h$, $k$,
respectively, then $(a,b,c)_\lambda\in\aA_{ghk}$. By
multilinearity, $\sS$ is graded. Intersecting with the original
graded subalgebra proves the assertion.
\end{proof}

\begin{lemma}\label{lem:idempotents_lift}
Let $\lambda\in\F\setminus\{-1\}$, with $\lambda=0$ allowed. Assume
that $\UTBV$ has a closed
$\Omega_\lambda$-grading for which $J(\UTBV)$ is graded. Every
homogeneous associative idempotent of $\UTBV/J(\UTBV)$ lifts to an
associative idempotent in $(\UTBV)_1$.
\end{lemma}

\begin{proof}
Set $\aA=\UTBV$ and $\iJ=J(\aA)$. The zero idempotent lifts to zero.
A nonzero homogeneous associative idempotent $\bar\ff$ has degree
$1$, since $\bar\ff\circ_\lambda\bar\ff=(1+\lambda)\bar\ff\neq0$.
Choose $x_0\in\aA_1$ mapping to $\bar\ff$, and define
\[
x_{n+1}=3x_n^2-2x_n^3
\quad\text{and}\quad
a_n=x_n^2-x_n.
\]
The identities
\[
x_{n+1}-x_n=-(2x_n-1)a_n
\quad\text{and}\quad
a_{n+1}=a_n^2(4a_n-3)
\]
give $a_n\in a_0^{2^n}\F[x_0]$. By
\Cref{lem:radicals:finite-algebras}, for each finite-dimensional
$W\subseteq\V$ there is an
$N$ such that $a_0^N(W)=0$. All the polynomials above commute with
$a_0$, so $a_n|_W=0$ and $x_{n+1}|_W=x_n|_W$ once $2^n\geq N$. The sequence
therefore converges in the complete algebra $\aA$ to an associative
idempotent $\ff$, with the same diagonal as $x_0$.

For a neutral element $x$, the identity
$x^{\circ_\lambda k}=(1+\lambda)^{k-1}x^k$ shows that every positive
associative power of $x$ is neutral. Thus $x_n\in\aA_1$ for every
$n$, and closedness gives $\ff\in\aA_1$. Its image is $\bar\ff$.
\end{proof}

\begin{corollary}\label{rem:primitive:idempotents:lift}
Under the hypotheses of \Cref{lem:idempotents_lift}, let
$\ee\in\UTBV$ be a neutral associative idempotent and let $\bar\ee$ be
its image in $\bar\aA=\UTBV/J(\UTBV)$. Then $\ee$ is primitive in
$(\UTBV)_1$ if and only if $\bar\ee$ is primitive in $\bar\aA_1$.
\end{corollary}
\begin{proof}
The corner $\ee\UTBV\ee$ is the fixed set of the continuous map
$x\mapsto\ee x\ee$, hence closed. By \Cref{rem:grad:inheritance}, its
homogeneous components are
$\ee(\UTBV)_g\ee=(\UTBV)_g\cap\ee\UTBV\ee$, which are closed, and the
corner meets $J(\UTBV)$ in $\ee J(\UTBV)\ee$. Hence the construction in the proof of
\Cref{lem:idempotents_lift} applies in this corner. If $\bar\ee$ is
not primitive, let $\bar\pp$ be a proper nonzero summand of it. A lift of $\bar\pp$ in $(\UTBV)_1$, multiplied on both sides by
$\ee$, is a neutral element of the corner mapping to $\bar\pp$, so
$\bar\pp$ lifts to an idempotent $\pp$ in the neutral component of
the corner. Then $\ee=\pp+(\ee-\pp)$ is a decomposition into nonzero
orthogonal idempotents in $(\UTBV)_1$. Conversely, any orthogonal
decomposition of $\ee$ in $(\UTBV)_1$ descends to a nontrivial
decomposition of its image, since the locally nilpotent ideal
$J(\UTBV)$ contains no nonzero idempotent.
\end{proof}

For $x\in\F^{\sI}$, write
$\Supp(x)=\{i\in\sI\mid x(i)\neq0\}$.

\begin{lemma}\label{lem:prim:idemp:comut}
Let $\sI$ be a set, and let $\aB=\F^{\sI}$ have an associative
$\G$-grading. For every $v\in\sI$ there is an idempotent
$\ff\in\aB_1$, primitive in $\aB_1$, with $\ff(v)\neq0$. Moreover, $\ff\aB\cong\F \H$ for a finite abelian
subgroup $\H\subseteq \G$ whose order is not divisible by
$\chr\F$. Every primitive idempotent of $\aB_1$ has finite support.
\end{lemma}

\begin{proof}
Write the coordinate idempotent $\ee_v$ as
$\ee_v=\sum_{h\in\H}\gamma_h$, where $\H\subseteq\G$ is finite, each
$0\neq\gamma_h\in\aB_h$, and put
$\aD=\Span\{\gamma_h\mid h\in \H\}$. For $b\in\aB_g$, compare
homogeneous components in
\[
\sum_{h\in \H}b\gamma_h=b\ee_v=b(v)\sum_{h\in \H}\gamma_h.
\]
If $b(v)=0$, then $b\aD=0$. If $b(v)\neq0$, the right side has
$|\H|$ nonzero homogeneous components. The left side has at most
$|\H|$ such components, of distinct degrees $gh$. Hence
$g\H=\H$ and $b\gamma_h=b(v)\gamma_{gh}$ for every $h\in \H$.

Choose $k\in \H$ with $\gamma_k(v)\neq0$. If $\gamma_h(v)=0$,
the first case gives $\gamma_h\gamma_k=0$, whereas the second case
gives $\gamma_k\gamma_h\neq0$, a contradiction. Thus every
$\gamma_h(v)$ is nonzero. It follows that $\H$ is a finite subgroup,
and that it is abelian because every product
$\gamma_h\gamma_k$ is nonzero. Applying the identity above to
$b=\gamma_h$ and $\gamma_1$ gives
$\gamma_h\gamma_1=\gamma_h(v)\gamma_h$. Evaluation at $v$ shows
that all $\gamma_h(v)$ equal $c=\gamma_1(v)$, and evaluation of
$\ee_v=\sum_h\gamma_h$ gives $|\H|c=1$.

The element $\ff=c^{-1}\gamma_1$ is an idempotent and the identity
of $\aD$. The component comparison also gives $b\ff\in\aD$ for
every homogeneous $b\in\aB$, so $\ff\aB=\aD$. The map
$h\mapsto c^{-1}\gamma_h$ identifies the canonically graded group
algebra $\F \H$ with $\aD$. Since
$\ff\aB\cong\F^{\Supp(\ff)}$ has dimension $|\H|$, the support of
$\ff$ is finite. Finally, $(\ff\aB)_1=\F\ff$, so $\ff$ is
primitive in $\aB_1$.

If $\pp\in\aB_1$ is any primitive idempotent, apply the construction
to $\pp\aB\cong\F^{\Supp(\pp)}$. It yields a nonzero neutral
idempotent $\qq\in\pp\aB$ of finite support. Primitivity forces
$\qq=\pp$, proving the last assertion.
\end{proof}

If $\G$ is torsion-free, then the grading of $\F^{\sI}$ is trivial.
Indeed, the preceding proof gives $\H=\{1\}$, so every coordinate
idempotent is neutral. For $b\in\aB_g$ with $g\neq1$, the identity
$b\ee_v=b(v)\ee_v$ then implies $b(v)=0$ for every $v$.

\begin{lemma}\label{lem:prim:idemp:existence}
Assume that either $\lambda\in\F\setminus\{0,\pm1\}$, or
$\lambda=1$ and $\chr\F\neq2$. Let $\UTBV$ have a closed
$\Omega_\lambda$-grading. For every $v\in\B$ there is an associative
idempotent $\ee\in(\UTBV)_1$, primitive in $(\UTBV)_1$, such that
$v\in\Supp(\ee)$ and $\Supp(\ee)$ is finite.
\end{lemma}

\begin{proof}
Set $\aA=\UTBV$ and $\iJ=J(\aA)$. By
\Cref{lem:jacobson:radical}, the quotient
$\aB=\aA/\iJ\cong\F^\B$ inherits an $\Omega_\lambda$-grading.
Since its product is commutative and $1+\lambda\neq0$, this is an
associative grading as well. By \Cref{lem:prim:idemp:comut}, there
is an idempotent $\bar\ee\in\aB_1$, primitive in $\aB_1$, with
finite support containing $v$. By
\Cref{lem:idempotents_lift,rem:primitive:idempotents:lift}, it lifts
to an idempotent $\ee\in\aA_1$ primitive in $\aA_1$. Since the
quotient map records diagonal coordinates,
$\Supp(\ee)=\Supp(\bar\ee)$.
\end{proof}

\begin{lemma}\label{lem:primitive:support:ut}
Assume that $\UTBV$ has a closed $\Omega_\lambda$-grading, and let
$\ee\in(\UTBV)_1$ be an associative idempotent, primitive in
$(\UTBV)_1$, with finite support.
\begin{enumerate}
\item\label{lem:primitive:support:ut:nonjordan}
If $\lambda\in\F\setminus\{0,\pm1\}$, then
$|\Supp(\ee)|=1$.
\item\label{lem:primitive:support:ut:jordan}
If $\lambda=1$ and $\chr\F\neq2$, then
$|\Supp(\ee)|\in\{1,2\}$. If $\B$ is infinite, the support is
a singleton.
\end{enumerate}
In each case where the support is a singleton, for every $v\in\B$
there is a primitive associative idempotent
$\ee_v\in(\UTBV)_1$ with $\Supp(\ee_v)=\{v\}$.
\end{lemma}

\begin{proof}
The finite support bounds follow from
\Cref{lem:primitive:support:finite:corner}. Suppose that
$\lambda=1$, $\chr\F\neq2$, and $\B$ is infinite. A finite family
$\ee_1,\ldots,\ee_n$ of pairwise orthogonal associative idempotents
primitive in $(\UTBV)_1$, each of finite support, can be enlarged.
Indeed, $\gg=1-\sum_{i=1}^n\ee_i$ is neutral by
\Cref{lem:identity:neutral} and has infinite support. By
\Cref{rem:grad:inheritance,prop:corner:iso:ut}, the
closed graded corner $\gg\UTBV\gg$ is topologically isomorphic to
the full triangular algebra on $\Supp(\gg)$, preserving diagonal
supports. By \Cref{lem:prim:idemp:existence}, this corner contains a
nonzero neutral associative idempotent $\ee_{n+1}$, primitive in
its neutral component, with finite support. It is orthogonal to
the preceding idempotents and is primitive in $(\UTBV)_1$: every
orthogonal decomposition of $\ee_{n+1}$ has its summands in
$\ee_{n+1}\UTBV\ee_{n+1}\subseteq\gg\UTBV\gg$.

Starting with $\ee_1=\ee$ gives an infinite orthogonal family to
which \Cref{lem:primitive:support:orthogonal} applies. Thus
$|\Supp(\ee)|=1$. The last assertion now follows from
\Cref{lem:prim:idemp:existence}.
\end{proof}

\begin{theorem}\label{thm:class}
Assume that either $\lambda\in\F\setminus\{0,\pm1\}$, or
$\lambda=1$, $\chr\F\neq2$, and $\B$ is infinite. If $\UTBV$ has a
closed $\Omega_\lambda$-grading, then $(\UTBV)_1$ contains a
\sset, and the grading is isomorphic to an elementary grading relative
to $\B$.
\end{theorem}

\begin{proof}
By \Cref{lem:primitive:support:ut}, for every $v\in\B$ there is a
neutral associative idempotent $\ee_v$ with $\Supp(\ee_v)=\{v\}$. Hence
\Cref{lem:finite:diagonalization} gives a \sset in $(\UTBV)_1$,
and \Cref{lem:gen:class:sset} gives the required isomorphism.
\end{proof}

\begin{corollary}
Under the hypotheses of \Cref{thm:class}, the grading has finite
support. If $\lambda=1$, its support generates an abelian subgroup
of $\G$.
\end{corollary}

\begin{proof}
By \Cref{thm:class}, the grading is isomorphic to an elementary
grading on $\UTBV$, which has finite support by
\Cref{lem:finite:grading}. In the Jordan case,
\Cref{prop:elementary-omega} applies to its restriction to
$\UTBVlim$, whose support is the same as that of the full algebra.
Graded isomorphisms preserve support.
\end{proof}

\section*{Declaration of generative AI and AI-assisted technologies in the manuscript preparation process}
During the preparation of this work the authors used AI-assisted tools
(OpenAI Codex and Anthropic Claude) to check proofs and to revise the
text. The authors reviewed and edited the content and take full
responsibility for the content of the publication.

\end{document}